\documentclass[12pt]{article}

\usepackage{lineno}
\nolinenumbers

\usepackage[utf8]{inputenc}
\usepackage[english]{babel}
\usepackage[table]{xcolor}
\usepackage[fleqn]{amsmath}
\usepackage{nccmath}
\usepackage{booktabs}
\usepackage{caption}
\usepackage{float}
\usepackage{sidecap}

\usepackage[backend=biber,style=alphabetic,sorting=ynt]{biblatex}

\usepackage{answers}
\usepackage{setspace}
\usepackage{graphicx}
\usepackage{enumitem}
\usepackage{multicol}
\usepackage{mathrsfs}
\usepackage[margin=1in]{geometry} 
\usepackage{amsmath,amsthm,amssymb}

\begin{document}

\title{An Asymptotic Form of the Generating Function \(\prod\limits_{k=1}^\infty \Big(1+\frac{x^k}{k}\Big)\)}
\author{Andreas B. G. Blobel\\ 
andreas.blobel@kabelmail.de} 
\maketitle

\begin{abstract}
\noindent It is shown that the sequence of rational numbers \(r(k)\) generated by the ordinary generating function\(\prod\limits_{k=1}^\infty \Big(1+\frac{x^k}{k}\Big)\) converges to a limit \(C > 0\). \(C\) can be expressed as \(C = \exp{\Big(- \sum\limits_{k = 2}^\infty \frac{(-1)^k}{k}\ \zeta(k) \Big)}\) where \(\zeta()\) denotes the Riemann zeta function.
\end{abstract}
\vspace{10mm}

\noindent The ordinary generating function (OGF)
\begin{ceqn}
\begin{subequations}
\begin{equation}\label{eq:PartFrac}
R(x)\ :=\ \prod\limits_{k=1}^\infty \Big(1+\frac{x^k}{k}\Big)\ =\ \sum\limits_{k=0}^\infty r(k)\ x^k  
\end{equation}
is closely related to the well known OGF
\begin{equation}\label{eq:PartCnt}
Q(x)\ :=\ \prod\limits_{k=1}^\infty \Big(1 + x^k\Big)\ =\ \sum\limits_{k=0}^\infty q(k)\ x^k
\end{equation}
\end{subequations}
\(Q(x)\) generates the sequence of counters for the number of integer partitions with distinct parts \cite{IntParts}. \(q(k)\) is equal to the number of partitions of \(k\) into distinct parts for each \(k \ge 0\) \cite{A009}.
\bigskip

\noindent A partition with distinct parts of integer \(k\) can be regarded as a finite set \(S\) of (distinct) positive integers \(i \ge 1\) whose sum equals \(k\).
Let \(\mathcal{P}(k)\) denote the set of all such partitions of \(k\) and let \(S \in \mathcal{P}(k)\). We then have

\begin{align}
    \sum\limits_{i \in S} i = k
\end{align}
With each partition $S \in \mathcal{P}(k)$ we can associate the inverse of the product of its (distinct) elements
\begin{align}\label{eq:IP}
    \textsc{ip}(S) := \frac{1}{\prod\limits_{i \in S} i}
\end{align}
\end{ceqn}

\noindent With this in mind $r(k)$ can be written as
\begin{subequations}
\begin{align}
    r(k) &= 1 &:\ r = 0\label{eq:SumIPa}\\
    r(k) &= \sum\limits_{S \in \mathcal{P}(k)} \textsc{ip}(S) = \displaystyle\sum\limits_{S \in \mathcal{P}(k)} \frac{1}{\prod\limits_{i \in S} i} &:\ r \ge 1\label{eq:SumIPb}
\end{align}
\end{subequations}
In other words, $r(k)$ is equal to the sum over all partitions \(S \in \mathcal{P}(k)\) of the reciprocal of the product of the elements of $S$ .

\vspace{10mm}
\noindent How does the sequence $r(k)$ given in \eqref{eq:SumIPa} and \eqref{eq:SumIPb} behave? Does it converge to some limit $C>0$? Taking the logarithm of \eqref{eq:PartFrac}, applying the Mercator series expansion \cite{Mercator}, and summing up columns first gives
\begin{fleqn} 
\begin{equation}\label{eq:LogR}
\begin{alignedat}{5}
&\ln{R(x)} = \sum_{k\ge1} \ln{\Big(1+\tfrac{x^k}{k}\Big)}
&&= x &&- \tfrac{1}{2}\ x^{2}
&&+ \tfrac{1}{3}\ x^{3} &&- \cdots\\[7pt]
& &&+ \tfrac{x^2}{2}
&&- \tfrac{1}{2} \Big[\tfrac{x^2}{2}\Big]^{2}
&&+ \tfrac{1}{3} \Big[\tfrac{x^2}{2}\Big]^{3} &&- \cdots\\[12pt]
& &&+ \tfrac{x^3}{3}
&&- \tfrac{1}{2} \Big[\tfrac{x^3}{3}\Big]^{2}
&&+ \tfrac{1}{3} \Big[\tfrac{x^3}{3}\Big]^{3} &&- \cdots\\[7pt]
& &&\ \vdots\\[7pt]
& &&= \text{Li}_1(x)
&&- \tfrac{1}{2}\ \text{Li}_2(x^2)
&&+ \tfrac{1}{3}\ \text{Li}_3(x^3) &&- \cdots\\[7pt]
\end{alignedat}
\end{equation}
\end{fleqn}
Here \(\text{Li}_s(x)\) denotes the so-called polylogarithm \cite{Polylogarithm}, a Dirichlet type series \cite{Dirichlet_series}.
\\ \\ \noindent We are looking for an asymptotic relation of the form
\begin{ceqn}
\begin{equation}\label{eq:rightarrow}
    R(x)\ \xrightarrow{x \to 1-}\ \frac{C}{1-x}
\end{equation}
for some constant \(C > 0\). This is equivalent to the existence of the limit
\begin{equation}\label{eq:Climes}
    C\ =\ \lim_{x \to 1-} (1-x)R(x)
\end{equation}
Taking the logarithm of \eqref{eq:Climes} gives
\begin{equation}\label{eq:logClimes}
    \ln{C}\ =\ \lim_{x \to 1-} \big(\ln{(1-x)} + \ln{R(x)}\big)
\end{equation}
If we insert \eqref{eq:LogR}, observe the identity \cite{PolylogarithmParticularValues}
\begin{equation}
    \text{Li}_1(x) = -\ln(1-x)
\end{equation}
\end{ceqn}
and finally set \(x = 1\), we arrive at the condition
\begin{equation}\label{eq:lnCzeta}
\begin{alignedat}{5}
    &\ln{C}\ &&=\ -\ \tfrac{1}{2}\ \text{Li}_2(1)\ &&+\ \tfrac{1}{3}\ \text{Li}_3(1)\ &&
    -\ \tfrac{1}{4}\ \text{Li}_4(1)
    &&+ \cdots\\[10pt]
    & &&=\ -\ \tfrac{1}{2}\ \zeta(2)\ &&+\ \tfrac{1}{3}\ \zeta(3)\ &&-\ \tfrac{1}{4}\ \zeta(4)\ &&+ \cdots
\end{alignedat}
\end{equation}
where \(\zeta(s)\) denotes the Riemann Zeta function \cite{Riemann_zeta_function}. We therefore have
\begin{ceqn}
\begin{equation}\label{eq:Czeta}
    C = \exp{\Bigg(- \sum_{k = 2}^\infty \frac{(-1)^k}{k}\ \zeta(k) \Bigg)}
\end{equation}
\end{ceqn}
\\ \\ \noindent We observe that \(\zeta(k)\) converges rapidly towards 1 \cite{Particular_values_of_the_Riemann_zeta_function}:

\begin{table}[H] 
\centering
\begin{tabular}{c|cc}
\toprule
 \(k\) & \multicolumn{2}{c}{\(\zeta(k) - 1\)}\\[7pt]
\midrule
2 & \(\frac{\pi^2}{6} - 1\) & 0.644934\\[5pt]
3 & - & 0.202057\\[5pt]
4 & \(\frac{\pi^4}{90} - 1\) & 0.082323\\[5pt]
5 & - & 0.036928\\[5pt]
6 & \(\frac{\pi^6}{945} - 1\) & 0.017343\\[5pt]
7 & - & 0.008349\\[5pt]
8 & \(\frac{\pi^8}{9450} - 1\) & 0.004077\\[5pt]
9 & - & 0.002008\\[5pt]
10 & \(\frac{\pi^{10}}{93555} - 1\) & 0.000995\\[5pt]
11 & - & 0.000494\\[5pt]
\bottomrule
\end{tabular}
\caption*{\(\zeta(k) \xrightarrow{k \to \infty} 1\)}
\label{table:field_zeta}
\end{table}

\vspace{8mm}
\noindent This motivates the decomposition of \eqref{eq:lnCzeta}
\begin{fleqn}
\begin{equation}\label{eq:lnCzetaDecomp}
\begin{alignedat}{5}
    \ln{C} &=\ &&-\ \tfrac{1}{2}\ \zeta(2)\ &&+\ \tfrac{1}{3}\ \zeta(3)
    &&-\ \tfrac{1}{4}\ \zeta(4)\ &&+\ \cdots\\[10pt]
    &=\ &&-\ \tfrac{1}{2}\ \Big[\zeta(2) - 1\Big]\ &&+\ \tfrac{1}{3}\ \Big[\zeta(3) - 1\Big]\ &&-\ \tfrac{1}{4}\ \Big[\zeta(4) - 1\Big]
    &&+\ \cdots\\[10pt]
    & &&-\ \tfrac{1}{2}\ &&+\ \tfrac{1}{3} &&-\ \tfrac{1}{4}
    &&+\ \cdots\\[10pt]
    &=\ &&-\ \Delta\ + \ln{2} - 1
\end{alignedat}
\end{equation}
where \(\Delta\) is defined as
\begin{equation}\label{eq:Delta}
\begin{alignedat}{5}
    \Delta &:=\ &&+\ \tfrac{1}{2}\ \Big[\zeta(2) - 1\Big]\ &&-\ \tfrac{1}{3}\ \Big[\zeta(3) - 1\Big]\ &&+\ \tfrac{1}{4}\ \Big[\zeta(4) - 1\Big]
    &&-\ \cdots
\end{alignedat}
\end{equation}
\end{fleqn}
\vspace{5mm}

\noindent We therefore have from \eqref{eq:lnCzetaDecomp}
\begin{ceqn}
\begin{equation}\label{eq:CzetaDelta}
    C = \frac{2}{e^{1+\Delta}}
\end{equation}
\end{ceqn}
\vspace{5mm}

\noindent From \eqref{eq:Delta} we derive the sequence of corrections \(\Delta_m\) as follows
\begin{ceqn}
\begin{align}\label{eq:Delta_m}
    \Delta_m &= \left\{
  \begin{array}{lr}
    0\ &:\ m = 1\\[10pt]
    \sum\limits_{k = 2}^{m} \frac{(-1)^k}{k} \big(\zeta(k) - 1\big)\ &:\ m \ge 2
  \end{array}\right.
\end{align}
\end{ceqn}
This creates the sequence
\begin{ceqn}
\begin{equation}
C_m = \frac{2}{\exp{(1+\Delta_m)}}\quad: m \ge 1
\end{equation}
\end{ceqn}
of approximations of \(C\) whose first elements are listed in table \ref{table:field_C}.

\begin{table}[H]
\centering
\begin{tabular}{rl|l}
\toprule
 \(m\) & \(\Delta_m\) & \(\frac{2}{\exp{(1+\Delta_m)}}\)\\[7pt]
\midrule
 1 & 0.0       & 0.7357589\\[5pt]
 2 & 0.3224670 & 0.5329542\\[5pt]
 3 & 0.2551147 & 0.5700863\\[5pt]
 4 & 0.2756955 & 0.5584734\\[5pt]
 5 & 0.2683100 & 0.5626133\\[5pt]
 6 & 0.2712005 & 0.5609894\\[5pt]
 7 & 0.2700078 & 0.5616589\\[5pt]
 8 & 0.2705174 & 0.5613727\\[5pt]
 9 & 0.2702943 & 0.5614980\\[5pt]
10 & 0.2703937 & 0.5614421\\[5pt]
11 & 0.2703488 & 0.5614674\\[5pt]
12 & 0.2703693 & 0.5614559\\[5pt]
13 & 0.2703599 & 0.5614612\\[5pt]
\bottomrule
\end{tabular}
\caption{Approximation of \(C\)}
\label{table:field_C}
\end{table}
\vspace{10mm}

\textbf{Useful recurrence relations for computation}
\vspace{5mm}

\noindent For \(n>0\) we define the finite products
\begin{ceqn}
\begin{subequations}
\begin{align}
    R_n(x)\ &:=\ \prod_{k=1}^n \Big(1+\frac{x^k}{k}\Big)
    \ =\ \sum\limits_{k=0}^\infty\ r_n(k)\ x^k\label{eq:FinitPartFrac}\\[7pt]
    Q_n(x)\ &:=\ \prod_{k=1}^n (1+x^k)
    \ =\ \sum\limits_{k=0}^\infty\ q_n(k)\ x^k\label{eq:FinitPartCnt}
\end{align}
\end{subequations}
\end{ceqn}

\noindent The integer numbers \(q_n(k)\) in \eqref{eq:FinitPartCnt} count the number of partitions of \(k\) with distinct parts \textit{where no part exceeds} \(n\).
The coefficients \(q_n(k)\) clearly have 3 basic properties:
\begin{subequations}
\begin{align}
    q_n(k)\ &=\ q(k) &&\text{if}\quad k\ <=\ n\label{eq:QbasicA}\\[7pt]
    q_n(k)\ &=\ 0 &&\text{if}\quad
    k\ >\ \tfrac{n(n+1)}{2}\label{eq:QbasicB}\\[7pt]
    \sum_{k\ge0}\ q_n(k)\ &=\ 2^n\label{eq:QbasicC}
\end{align}
\end{subequations}
where \eqref{eq:QbasicC} follows from evaluation of \(Q_n(1)\).
The \(q_n(k)\) obey the recurrence relations
\begin{subequations}
\begin{align}
    q_0(k) &= \left\{
  \begin{array}{lr}
    1\ :\ k = 0\\
    0\ :\ k \ge 1
  \end{array}\right.\label{eq:QrecurA}\\[7pt]
    q_n(k) &=\ q_{n-1}(k)\quad:\ 0 \le k < n\label{eq:QrecurB}\\[7pt]
    q_n(k) &=\ q_{n-1}(k-n)\ +\ q_{n-2}(k-n+1)\ +\ q_{n-3}(k-n+2)
    \ +\ \cdots\nonumber\\[7pt]
    &\ +\ q_{1}(k-2)\ +\ q_{0}(k-1)\quad:\ k \ge n > 0\label{eq:QrecurC}
\end{align}
\end{subequations}
Initial values are prescribed in row \(n = 0\) \eqref{eq:QrecurA}. The values in any subsequent row \(n \ge 1\) are determined by values in previous rows \(m < n\).
\vspace{6mm}

\begin{table}[H]
\centering
\begin{tabular}{r|ccccccccccccccccc}
  & 0 & 1 & 2 & 3 & 4 & 5 & 6 & 7 & 8 & 9 & 10 & 11 & 12 & 13 & 14 & 15 & 16\\
\hline
0 & \cellcolor{orange}1 & 0 & 0 & 0 & 0 & 0 & 0 & 0 & 0 & 0 & 0 & 0 & 0 & 0 & 0 & 0 & 0\\[4pt]
1 & 1 & \cellcolor{orange}1 & 0 & 0 & 0 & 0 & 0 & 0 & 0 & 0 & 0 & 0 & 0 & 0 & 0 & 0 & 0\\[4pt]
2 & 1 & 1 & 1 & \cellcolor{orange}1 & 0 & 0 & 0 & 0 & 0 & 0 & 0 & 0 & 0 & 0 & 0 & 0 & 0\\[4pt]
3 & 1 & 1 & 1 & 2 & 1 & 1 & \cellcolor{orange}1 & 0 & 0 & 0 & 0 & 0 & 0 & 0 & 0 & 0 & 0\\[4pt]
4 & 1 & 1 & 1 & 2 & 2 & 2 & 2 & 2 & 1 & 1 & \cellcolor{orange}1 & 0 & 0 & 0 & 0 & 0 & 0\\[4pt]
5 & 1 & 1 & 1 & 2 & 2 & 3 & 3 & 3 & 3 & 3 & 3 & 2 & 2 & 1 & 1 & \cellcolor{orange}1 & 0\\[4pt]
\end{tabular}
\caption{Upper left section of the \(q_n(k)\) field \([\ 0 \le n \le 5\ ,\ 0 \le k \le 16\ ]\)}
\label{table:field_q}
\end{table}

\vspace{10mm}
\noindent Analogous properties and relations hold for the \textit{rational} numbers $r_n(k)$ in \eqref{eq:FinitPartFrac}:
\begin{subequations}
\begin{align}
    r_n(k)\ &=\ r(k) &&\text{if}\quad k\ <=\ n\label{eq:RbasicA}\\[7pt]
    r_n(k)\ &=\ 0 &&\text{if}\quad
    k\ >\ \tfrac{n(n+1)}{2}\label{eq:RbasicB}\\[7pt]
    \sum_{k\ge0}\ r_n(k)\ &=\ n+1\label{eq:RbasicC}
\end{align}
\end{subequations}
\vspace{3mm}

\begin{subequations}
\begin{align}
    r_0(k) &= \left\{
  \begin{array}{lr}
    1\ :\ k = 0\\
    0\ :\ k \ge 1
  \end{array}\right.\label{eq:RrecurA}\\[7pt]
    r_n(k) &=\ r_{n-1}(k)\quad:\ 0 \le k < n\label{eq:RrecurB}\\[7pt]
    r_n(k) &=\ \tfrac{1}{n}\ r_{n-1}(k-n)\ +\ \tfrac{1}{n-1}\ r_{n-2}(k-n+1)
    \ +\ \tfrac{1}{n-2}\ r_{n-3}(k-n+2) \ +\ \cdots\nonumber\\[7pt]
    &\ +\ \tfrac{1}{2}\ r_{1}(k-2)\ +\ \tfrac{1}{1}\ r_{0}(k-1)\quad:\ k \ge n > 0\label{eq:RrecurC}
\end{align}
\end{subequations}
\vspace{5mm}

\begin{table}[H]
\centering
\begin{tabular}{r|ccccccccccccccccc}
  & 0 & 1 & 2 & 3 & 4 & 5 & 6 & 7 & 8 & 9 & 10 & 11 & 12 & 13 & 14 & 15 & 16\\
\hline
0 & \cellcolor{orange}1 & 0 & 0 & 0 & 0 & 0 & 0 & 0 & 0 & 0 & 0 & 0 & 0 & 0 & 0 & 0 & 0\\[6pt]
1 & 1 & \cellcolor{orange}1 & 0 & 0 & 0 & 0 & 0 & 0 & 0 & 0 & 0 & 0 & 0 & 0 & 0 & 0 & 0\\[6pt]
2 & 1 & 1 & \(\frac{1}{2}\) & \cellcolor{orange}\(\frac{1}{2}\) & 0 & 0 & 0 & 0 & 0 & 0 & 0 & 0 & 0 & 0 & 0 & 0 & 0\\[6pt]
3 & 1 & 1 & \(\frac{1}{2}\) & 
\(\frac{5}{6}\) & \(\frac{1}{3}\) & \(\frac{1}{6}\) & \cellcolor{orange}\(\frac{1}{6}\) & 0 & 0 & 0 & 0 & 0 & 0 & 0 & 0 & 0 & 0\\[6pt]
4 & 1 & 1 & \(\frac{1}{2}\) & \(\frac{5}{6}\) & \(\frac{7}{12}\) & \(\frac{5}{12}\) & \(\frac{7}{24}\) & \(\frac{5}{24}\) & \(\frac{1}{12}\) & \(\frac{1}{24}\) & \cellcolor{orange}\(\frac{1}{24}\) & 0 & 0 & 0 & 0 & 0 & 0\\[6pt]
5 & 1 & 1 & \(\frac{1}{2}\) & \(\frac{5}{6}\) & \(\frac{7}{12}\) & \(\frac{37}{60}\) & \(\frac{59}{120}\) & \(\frac{37}{120}\) & \(\frac{1}{4}\) & \(\frac{19}{120}\) & \(\frac{1}{8}\) & \(\frac{7}{120}\) & \(\frac{1}{24}\) & \(\frac{1}{60}\) & \(\frac{1}{120}\) & \cellcolor{orange}\(\frac{1}{120}\) & 0\\[6pt]
\end{tabular}
\caption{Upper left section of the \(r_n(k)\) field \([\ 0 \le n \le 5\ ,\ 0 \le k \le 16\ ]\)}
\label{table:field_r}
\end{table}
\vspace{10mm}

\begin{figure}[H]
\caption{Some computed instances of \(r(k)\)}
\begin{center}
      \pdfimage{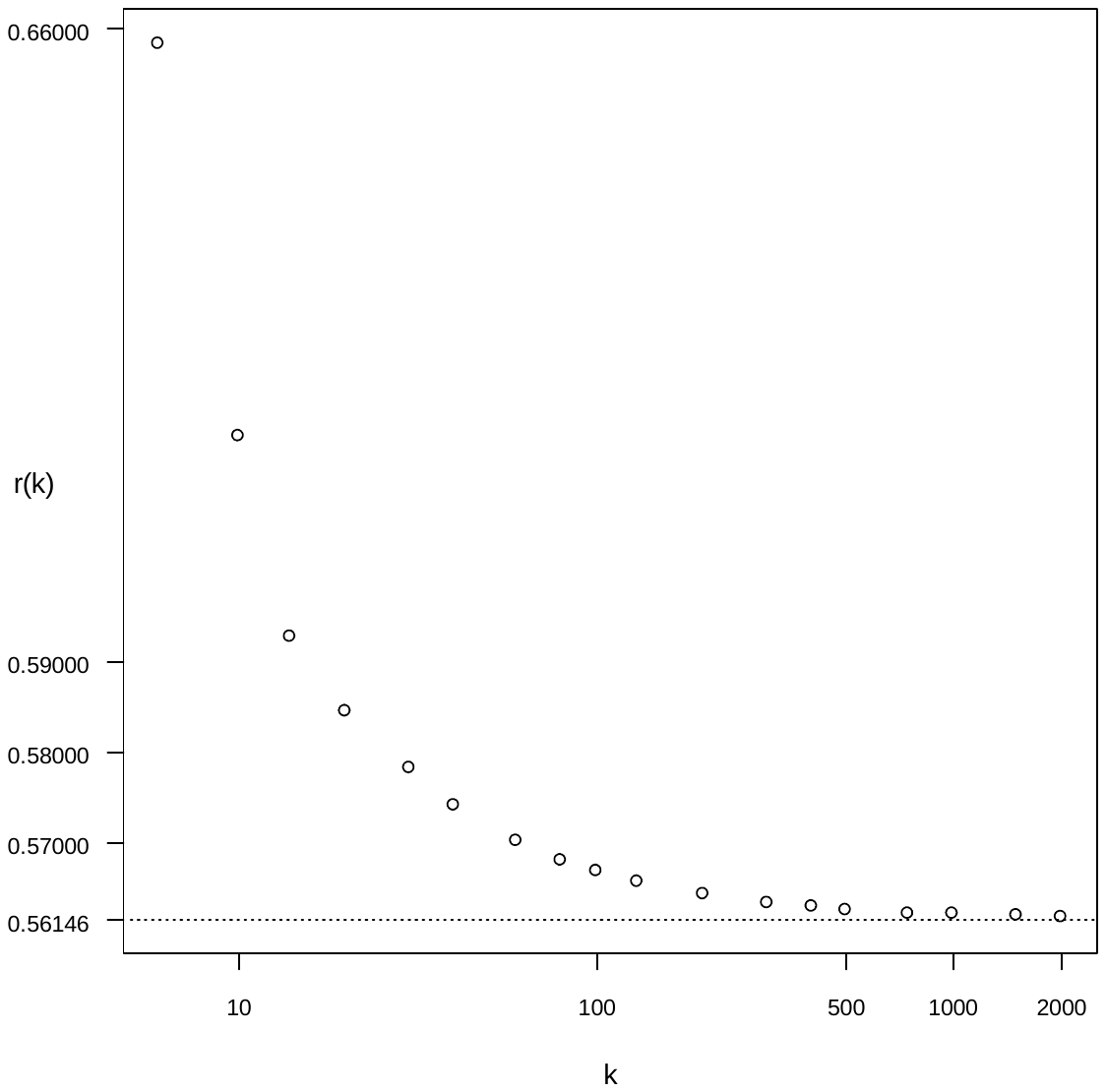}
\end{center}
\label{ComputedValues}
\end{figure}

\noindent Figure \ref{ComputedValues} assembles some instances of $r(k)$ which have been computed on the R platform for statistical computing \cite{R_Project} using recurrence relations \eqref{eq:RrecurA}, \eqref{eq:RrecurB}, and \eqref{eq:RrecurC}.
The plot shows that the \(r(k)\) approach the asymptotic value
\begin{ceqn}
\[C = 0.56146\dots\]
\end{ceqn}
from above as \(k\) increases. The constant \(C\) is determined by \eqref{eq:Czeta} and \eqref{eq:CzetaDelta} and is marked by a dashed horizontal line.

\textbf{Conclusion}
\vspace{5mm}

\noindent It has been shown that the function
\begin{ceqn}
\[f(x) = \frac{C}{1-x}\]
\end{ceqn}
is an asymptotic form of the generating function \eqref{eq:PartFrac} in the sense that the sequence of rational numbers \(r(k)\) generated by \eqref{eq:PartFrac} converges towards \(C > 0\) which is determined by \eqref{eq:Czeta} and \eqref{eq:CzetaDelta}.

\pagebreak

\printbibliography

\end{document}